\documentclass[12pt]{article}
\usepackage{amssymb}
\usepackage{mathrsfs}%有更多的数学符号，比如花体$\mathscr P$
\usepackage{pst-poly}     % From pstricks/contrib/pst-poly
\usepackage{cite}
\usepackage{graphicx}
\usepackage{pst-plot}  %坐标
\usepackage{tikz}

\newtheorem{thm}{Theorem}[section]
\newtheorem{lem}[thm]{Lemma}
\newtheorem{Def}[thm]{Definition}

\newtheorem{cor}[thm]{Corollary}

\newtheorem{rem}[thm]{Remark}

\newtheorem{prob}[thm]{Problem}

\newenvironment{pf}[1][Proof]{\noindent\textbf{#1.} }{\hfill\rule{1mm}{2mm}}

\makeatletter \@addtoreset{equation}{section} \makeatother

\begin{document}
\title{\bf On fault tolerance of $(n,k)$-star networks\thanks
{The work was supported by NNSF of China(10711233, 11571044, 11601041, 61673006),
Young talent fund from Hubei EDU(Q20151311) and Yangtze University(2015cqr23), NSF of Hubei(2014CFB248).}}

\author
{Xiang-Jun Li$^a$ \quad Yong-ni Guan$^a$\quad Zheng Yan$^a$ \quad
Jun-Ming Xu$^b$ \footnote{Corresponding author: xujm@ustc.edu.cn
(J.-M. Xu) }\\
\\
 {\small $^a$School of Information and Mathematics,}\\
 {\small Yangtze University, Jingzhou, Hubei, 434023, China}\\
 {\small $^b$School of Mathematical Sciences}\\
 {\small University of Science and Technology of China, Hefei, 230026, China}}
\date{}
\maketitle

%\begin{frontmatter}
\begin{abstract}

Fault tolerance of an $(n,k)$-star network is measured by its
$h$-super connectivity $\kappa_s^{(h)}$ or $h$-super
edge-connectivity $\lambda_s^{(h)}$. Li {\it et al.} [Appl. Math.
Comput. 248 (2014), 525-530; Math. Sci. Lett. 1 (2012), 133-138]
determined $\kappa_s^{(h)}$ and $\lambda_s^{(h)}$ for $0\leq h\leq
n-k$. This paper determines
$\kappa_s^{(h)}=\lambda_s^{(h)}=\frac{(h+1)!(n-h-1)}{(n-k)!}$ for
$n-k\leq h \leq n-2$.

\vskip6pt


\end{abstract}

\section{Introduction}
It is well known that interconnection networks play an important
role in parallel computing/communication systems. An interconnection
network can be modeled by a graph in which vertices correspond to
processors and edges correspond to communication links.

Let $G$ be a connected graph. A subset $T\subset V(G)$, if any, is
called an {\it $h$-vertex-cut}, if $G-T$ is disconnected and has the
minimum degree at least $h$. The {\it $h$-super connectivity}
$\kappa_s^{(h)}(G)$ of $G$ is defined as the minimum cardinality
over all $h$-vertex-cuts of $G$. Similarly, a subset $F\subset
E(G)$, if any, is called an {\it $h$-edge-cut}, if $G-F$ is
disconnected and has the minimum degree at least $h$. The {\it
$h$-super edge-connectivity} $\lambda_s^{(h)}(G)$ of $G$ is defined
as the minimum cardinality over all $h$-edge-cuts of $G$.

The $h$-super connectivity and $h$-super edge-connectivity are
important measure of fault tolerance of networks and have been
received considerable attention in the literature (see, for example,
\cite{lx12, lxu14, lx14} and references cited therein).

For the $n$-dimensional star graph $S_{n}$, Li and Xu~\cite{lx14}
proved that $\kappa_s^{(h)}(S_{n})=\lambda _s^{(h)}(S_{n})
=(h+1)!(n-h-1)$ for any $h$ with $0\leq h\leq n-2$. As a
generalization of $S_n$, the $(n,k)$-star graph $S_{n,k}$, where
$2\leq k\leq n-1$, when $0\leq h\leq n-k$, Li and Xu~\cite{lx12,
lxu14} determined that
 \begin{equation}\label{e1.1}
 \kappa _s^{(h)}(S_{n,k})=n+h(k-2)-1
 \end{equation}
and
  \begin{equation}\label{e1.2}
 \lambda_s^{(h)}(S_{n,k})= \left\{\begin{array}{ll}
 (n-h-1)(h+1)&\ {\rm for}\ h\leqslant \min\{k-2,\frac{n}{2}-1\},\\
 (n-k+1)(k-1)&\ {\rm otherwise}.
    \end{array}\right.
 \end{equation}
However, when $n-k+1 \leq h\leq n-2$, $\kappa_s^{(h)}(S_{n,k})$ and
$\lambda_s^{(h)}(S_{n,k})$ have not been determined yet. In this
paper, we prove that
 $$
 \kappa_s^{(h)}(S_{n,k})=\lambda_s^{(h)}(S_{n,k})=\frac{(h+1)!(n-h-1)}{(n-k)!}
 $$
by using an $(n-k)!$-split graph of $S_{n,k}$.

The rest of the paper is organized as follows. In Section 2, we give
definitions of a star graph $S_n$, a $(n,k)$-star graph $S_{n,k}$
and an $(n-k)!$-split graph of $S_{n,k}$, and some lemmas used in
our proofs. The proof of our main result is in Section 3.
Conclusions and problems are in Section 4.

\section{Definitions and lemmas}

For a given integer $n$ with $n\geq2$, let $I_n=\{1,2,\ldots,n\}$,
$I'_n=\{2,\ldots,n\}$. For $k\in I_n$, let $P(n,k)$ be the set of
$k$-arrangements on $I_n$, that is, $P(n,k)=\{ p_{1}p_{2}\ldots
p_{k}:\ p_{i}\in I_n, p_{i}\neq p_{j}, 1\leq i\neq j\leq k\}$.
$P(n,n)$ will be shorted as $P(n)$. Clearly,
$|P(n,k)|=\frac{n!}{(n-k)!}$. Usually, if $u=p_1p_2\dots p_k\in
P(n,k)$, we call $p_i$ the {\it $i$-digit} of $u$ for each $i\in
I_k$. For simplicity, we write $up_{k+1}\ldots p_{n}$ for
$p=p_1p_2\ldots p_k p_{k+1}\ldots p_{n}\in P(n)$, where $u$ is
called the prefix of $p$ and $p_{k+1}\ldots p_{n}$ is called the
suffix of $p$.

\begin{Def}\label{def2.1} \textnormal{(Akers and Krishnamurthy~\cite{ak89}, 1989)}
An $n$-dimensional star graph $S_{n}$ is a graph with vertex-set
$P(n)$, a vertex $p=p_{1}p_{2}\ldots p_{i}\ldots p_{n}$ being linked
a vertex $q$ if and only if $q=p_{i}p_{2}\ldots
p_{i-1}p_1p_{i+1}\ldots p_{n}$ for some $i\in I'_n$.
  \end{Def}

\begin{lem}\label{lem2.1}\textnormal{(Li and Xu~\cite{lx14}, 2014)}
$\kappa_s^{(h)}(S_{n})=\lambda_s^{(h)}(S_{n})=(h+1)!(n-h-1)$ for any
$h$ with $0 \leq h \leq n-2$.
\end{lem}

\begin{Def}\label{def2.2} \textnormal{(Chiang {\it et al.}~\cite{cc95}, 1995)}
An $(n,k)$-star graph $S_{n,k}$ is a graph with vertex-set $P(n,k)$,
a vertex $p=p_{1}p_{2}\ldots p_{i}\ldots p_{k}$ being linked a
vertex $q$ if and only if $q$ is

{\rm (a)}\ $p_{i}p_{2}\cdots p_{i-1}p_{1}p_{i+1}\cdots p_{k}$, where
$i\in I'_k$ (swap $p_{1}$ with $p_{i}$), or

{\rm (b)} $p'_{1} p_{2}p_{3}\cdots p_{k}$, where $p'_{1}\in I_n\setminus \{p_{i}:\
i \in I_k \}$ (replace $p_{1}$ by $p'_{1}$).
\end{Def}

The vertices of type $(a)$ are referred to as {\it swap-neighbors}
of the vertex $p$ and the edges between them are referred to as {\it
swap-edges} or {\it $i$-edges}. The vertices of type $(b)$ are
referred to as {\it unswap-neighbors} of the vertex $p$ and the
edges between them are referred to as {\it unswap-edges}. Clearly,
every vertex in $S_{n,k}$ has $(k-1)$ swap-neighbors and $(n-k)$
unswap-neighbors.

By definitions, it is clear that $S_{n,1}\cong K_n$, a complete
graph with $n$ vertices, and $S_{n, n-1}\cong S_n$.

\begin{Def}\
Let $G$ be a graph and $t$ be a positive integer. A {\it $t$-split}
graph $G^t$ of $G$ is a graph obtained from $G$ by replacing each
vertex $x$ by a set $V_x$ of $t$ independent vertices, and replacing
each edge $e=xy$ by a perfect matching $E_e$ between $V_x$ and
$V_y$.
\end{Def}

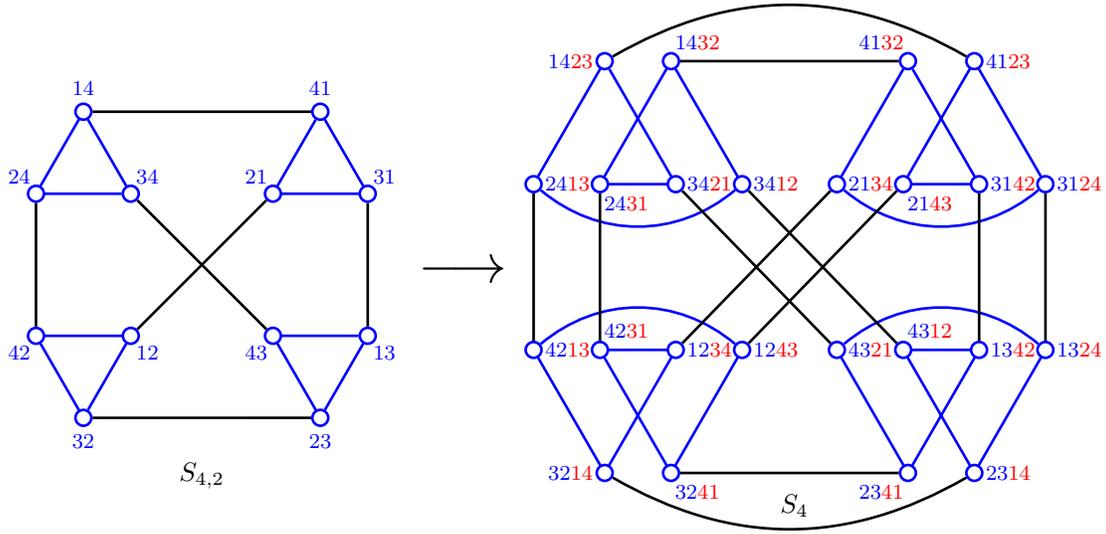
\begin{figure}[ht]
\begin{center}
\begin{tikzpicture}[scale=.63]
    \tikzstyle{every node}=[draw,circle,fill=white,minimum size=6pt,
                            inner sep=0pt]

    \draw[line width = 1pt,blue] (-2,-1.8) node (34) [label=45:\scriptsize$  3 4$] {}
        -- ++(120:2.0cm) node (14) [label=90:\scriptsize$ 14$] {}
        -- ++(240:2.0cm) node (24) [label=135:\scriptsize$24$] {}
        -- (34);
   \draw[line width = 1pt,blue] (14)
         ++(0:5.0cm) node (41) [label=90:\scriptsize$41$]    {}
        -- ++(-60:2.0cm) node (31) [label=45:\scriptsize$ 31$] {}
        -- ++(180:2.0cm) node (21) [label=135:\scriptsize$21$] {}
        -- (41);
   \draw[line width = 1pt,blue] (31)
         ++(-90:3.0cm) node (13) [label=-45:\scriptsize$13$]  {}
        -- ++(240:2cm) node (23) [label=-90:\scriptsize$23$]    {}
        -- ++(120:2.0cm) node (43) [label=-135:\scriptsize$43$] {}
        -- (13);
    \draw[line width = 1pt,blue] (23)
        ++(180:5.0cm) node (32) [label=-90:\scriptsize$32$]  {}
        -- ++(120:2.0cm) node (42) [label=-135:\scriptsize$42$] {}
        -- ++(0:2.0cm) node (12) [label=-45:\scriptsize$12$]    {}
        -- (32);
    % Add missing ``straight'' edges
    \draw[line width = 1pt] (42) -- (24);
    \draw[line width = 1pt] (43) -- (34);
    \draw[line width = 1pt] (12) -- (21);
    \draw[line width = 1pt] (14) -- (41);
    \draw[line width = 1pt] (13) -- (31);
    \draw[line width = 1pt] (23) -- (32);

   \node [line width = 1.0pt, white] at (-0.5, -7)
     [label=-90:  {\footnotesize $S_{4,2}$} ]{};

  %\the right

   \node [line width = 1.0pt, white] at (5, -2.3)
     [label=-90:  {\LARGE $\longrightarrow$} ]{};
   \node [line width = 1.0pt, white] at (12, -7.8)
     [label=-90:  {\footnotesize $S_{4}$} ]{};
    % First, draw the inner hexagon with a ``pin'' -- namely,
   \draw [line width = 1pt,blue](8,1) node (1423) [label=180:{\scriptsize 14{\red  23}}] {}
        -- ++(-120:3cm) node (2413) [label=2: {\scriptsize 24{\red13}}] {}
         ++(0:4.4cm) node (3412) [label=right: {\scriptsize 34{\red12}}] {}
        -- ++(120:3cm) node (1432) [label=10: {\scriptsize 14{\red32}}] {}
        -- ++(-120:3cm) node (2431) [label=-10: {\scriptsize 24{\red31}}] {}
        -- ++(0:1.6cm) node (3421) [label=right:
        {\scriptsize 34{\red21}}] {}
        -- (1423); % 1  shape is closed, we now connect it to an outer vertex:
    \draw [line width = 1pt,blue](2413) to [out=-40,in=220] (3412);

     \draw [line width = 1pt,blue](8,1) %node (2431) {}
         ++(0:1.4cm) %node (1432) {}
         ++(0:5.0cm) node (4132) [label=170: {\scriptsize 41{\red32}}] {}
       -- ++(-120:3cm) node (2134) [label=0: {\scriptsize 21{\red34}}] {}
         ++(0:4.4cm) node (3124) [label=right: {\scriptsize 31{\red24}}] {}
        -- ++(120:3cm) node (4123) [label=0: {\scriptsize 41{\red23}}] {}
        -- ++(-120:3cm) node (2143) [label=-10: {\scriptsize 21{\red43}}] {}
        -- ++(0:1.6cm) node (3142) [label=right: {\scriptsize 31{\red42}}] {}
        -- (4132); % 2  shape is closed, we now connect it to an outer vertex:
        \draw [line width = 1pt,blue](2134) to [out=-40,in=220] (3124);

        \draw [line width = 1pt,blue](8,1) %node (2431) {}
        ++(-120:3cm) ++(-90:3.5cm) node (4213) [label=right: {\scriptsize 42{\red13}}] {}
        -- ++ (-60:3cm) node (3214) [label=left: {\scriptsize 32{\red14}}] {}
         -- ++ (60:3cm) node (1234) [label=right: {\scriptsize 12{\red34}}] {}
         -- ++(180:1.6cm) node (4231) [label=10: {\scriptsize 42{\red31}}] {}
        -- ++ (-60:3cm) node (3241) [label=-10: {\scriptsize 32{\red41}}] {}
        -- ++ (60:3cm) node (1243) [label=right: {\scriptsize 12{\red43}}] {};
        \draw [line width = 1pt,blue](4213) to [out=40,in=140] (1243);
         % 3 shape is closed, we now connect it to an outer vertex:

     \draw [line width = 1pt,blue](8,1) %node (2431) {}
        ++(-120:3cm) ++(-90:3.5cm) ++(0:6.4cm)
         node (4321) [label=right: {\scriptsize 43{\red21}}] {}
        -- ++ (-60:3cm) node (2341) [label=-170: {\scriptsize 23{\red41}}] {}
         -- ++ (60:3cm) node (1342) [label=right: {\scriptsize 13{\red42}}] {}
         -- ++(180:1.6cm) node (4312) [label=10: {\scriptsize 43{\red12}}] {}
        -- ++ (-60:3cm) node (2314) [label=right: {\scriptsize 23{\red14}}] {}
        -- ++ (60:3cm) node (1324) [label=right: {\scriptsize 13{\red24}}] {};
        \draw [line width = 1pt,blue](4321) to [out=40,in=140] (1324);
         % 4 shape is closed, we now connect it to an outer vertex:

    \draw [line width = 1pt,black](2413) -- (4213);
    \draw [line width = 1pt,black](2431) -- (4231);
    \draw [line width = 1pt,black](3421) -- (4321);
    \draw [line width = 1pt,black](3412) -- (4312);
    \draw [line width = 1pt,black](2134) -- (1234);
    \draw [line width = 1pt,black](2143) -- (1243);
    \draw [line width = 1pt,black](3142) -- (1342);
    \draw [line width = 1pt,black](3124) -- (1324);

    \draw [line width = 1pt,black](1432) -- (4132);
    \draw [line width = 1pt,black](3241) -- (2341);

    \draw [line width = 1pt,black](1423) to [out=30,in=150] (4123);
    \draw [line width = 1pt,black](3214) to [out=-30,in=-150] (2314);

\end{tikzpicture}
\end{center}
\vskip-.8cm \caption{\label{f1}\footnotesize A $(4,2)$-star graph
$S_{4,2}$ and its $2$-split graph $S^2_{4,2}$, which is isomorphic
to a star $S_{4}$. }
\end{figure}

Fig.~\ref{f1} shows a $(4,2)$-star graph $S_{4,2}$ and its $2$-split
graph $S^2_{4,2}$, which is isomorphic to a star $S_{4}$.

\begin{lem}\label{lem2.5}
Let $G$ be a connected graph and $G^t$ be a $t$-split graph of $G$.
Then $\kappa_s^{(h)}(G^t)\leq t\,\kappa_s^{(h)}(G)$ and
$\lambda_s^{(h)}(G^t)\leq t\,\lambda_s^{(h)}(G)$.
\end{lem}

\begin{pf}
Assume that $T$ is a minimum $h$-vertex-cut and $F$ is a minimum
$h$-edge-cut in $G$. Then $\kappa_s^{(h)}(G)= |T|$  and
$\lambda_s^{(h)}(G)= |F|$. Let $T^t=\{V_u:\ u\in T \}$ and
$F^t=\{E_e:\ e\in F\}$. Then $|T^t|=t\,|T|$ and $|F^t|=t\,|F|$.

Since $G-T$ (resp. $G-F$) is disconnected, then $G^t-T^t$ (resp.
$G^t-F^t$) also is disconnected. Furthermore, it is easy to see that
$(G-T)^t=G^t-T^t$ (resp. $(G-F)^t=G^t-F^t$).

Because $T$ (resp. $F$) is an $h$-vertex-cut (resp. $h$-edge-cut) in
$G$, each vertex in $G-T$ (resp. $G-F$) has at least $h$ neighbors
in $G-T$ (resp. $G-F$), and so each vertex in $(G-T)^t$ (resp.
$(G-F)^t$) also has $h$ neighbors in $(G-T)^t$ (resp. $(G-F)^t$),
which implies that $T^t$ is an $h$-vertex-cut (resp. $F^t$ is an
$h$-edge-cut) in $G^t$. Thus, we have
 $$
 \kappa_s^{(h)}(G^t)\leq |T^t|=t\, |T|= t\,\kappa_s^{(h)}(G),
 $$
 $$
 \lambda_s^{(h)}(G^t)\leq |F^t|=t\, |F|= t\,\lambda_s^{(h)}(G)
 $$
as required.
\end{pf}

\begin{lem}\label{lem2.6}
For any $k$ with $2 \leq k\leq n-1$, there is an $(n-k)!$-split
graph of $S_{n,k}$ that is isomorphic to a star graph $S_n$.
\end{lem}

\begin{pf}
Define an $(n-k)!$-split graph $S^{(n-k)!}_{n,k}$ of $S_{n,k}$ as
follows.

For a vertex $u=p_1p_2\ldots p_k$ in $S_{n,k}$, it is replaced by
$(n-k)!$ vertices
$$
V_u=\{up_{k+1}\ldots p_{n}\in P(n):\  p_{k+i}\in I_n\setminus \{
p_1,\ldots, p_k\} \textnormal{\ for\ }  1\leq i\leq n-k\}.
$$

For an edge $uv$ in $S_{n,k}$, let $x=up_{k+1}\ldots p_{n}\in V_u$,
and define a matching $E_{uv}$ between $V_u$ and $V_v$ as follows.

If $uv$ is an $i$-edge in $S_{n,k}$ for some $i\in I'_k$, then
$v=p_ip_2\ldots p_{i-1}p_1p_{i+1}\ldots p_k$. Let $E_{uv}$ be the
set of edges that link two vertices $x\in V_u$ and $y\in V_v$ with
the same suffix.

If $uv$ is an unswap-edge in $S_{n,k}$, then $v=p_{k+j}
p_{2}p_{3}\cdots p_{k}$ for some $p_{k+j}\in I_n\setminus \{p_{i}:\
i \in I_k \}$. Let $E_{uv}$ be the set of edges that link two
vertices $x\in V_u$ and $y\in V_v$ with suffixes differing in
exactly the $(k+j)$-digit.

Clearly, $S^{(n-k)!}_{n,k}$ has vertex-set $P(n)$, a vertex $x$ is
adjacent to a vertex $y$ if and only if the label of $y$ can be
obtained from the label of $x$ by swapping the first digit and the
$i$-digit for some $i\in I'_n$. Therefore, by
Definition~\ref{def2.1}, $S^{(n-k)!}_{n,k}$ is a star graph $S_n$.
The Lemma follows.
\end{pf}

\section{Main results}

In this section, we present our main results, that is, we determine
the $h$-super connectivity and $h$-super edge connectivity of the
$(n,k)$-star graph $S_{n,k}$. Since $S_{n,1}\cong K_n$, for which
$\kappa^{(h)}_s$ and $\lambda^{(h)}_s$ do not exist for any $h$ with
$1\leq h\leq n-1$, we only consider the case of $k\geq 2$ in the
following discussion.

\begin{lem}\label{lem3.1}
For $2 \leq k\leq n-1$ and $n-k\leq h \leq n-2$,
 $$
 \lambda_s^{(h)}(S_{n,k})\geq \frac{(h+1)!(n-1-h)}{(n-k)!}\ {\rm and}\ \kappa_s^{(h)}(S_{n,k})\geq \frac{(h+1)!(n-1-h)}{(n-k)!}.
 $$
\end{lem}

\begin{pf}
For $2 \leq k \leq n-1$ and $n-k\leq h \leq n-2$, by
Lemma~\ref{lem2.5}, Lemma~\ref{lem2.6} and  Lemma~\ref{lem2.1}, we
immediately have that
 $$
  \begin{array}{rl}
 &\kappa_s^{(h)}(S_{n,k})(n-k)! \geq
 \kappa_s^{(h)}(S^{(n-k)!}_{n,k})=\kappa_s^{(h)}(S_{n})=
 (h+1)!(n-h-1)\\
 &\lambda_s^{(h)}(S_{n,k})(n-k)! \geq
 \lambda_s^{(h)}(S^{(n-k)!}_{n,k})=\lambda_s^{(h)}(S_{n})=
 (h+1)!(n-h-1)
 \end{array}
 $$
as required.
\end{pf}

\begin{lem}\label{lem3.2}
For $2 \leq k\leq n-1$ and $n-k\leq h \leq n-2$,
 $$
 \lambda_s^{(h)}(S_{n,k})\leq \frac{(h+1)!(n-1-h)}{(n-k)!}\ {\rm and}\
\kappa_s^{(h)}(S_{n,k})\leq \frac{(h+1)!(n-1-h)}{(n-k)!}.
 $$
\end{lem}

\begin{pf}
Since $n-k\leq h$, $n-1-h\leq k-1$. Let $X$ be the set of
$k$-arrangements on $I_n$ whose the last $(n-1-h)$ digits are
$12\cdots(n-1-h)$. Then $|X|=\frac{(h+1)!}{(n-k)!}$. Let $H$ be the
subgraph of $S_{n,k}$ induced by $X$. Since $n\geq k+1$,
$h+1-(n-k)\leq h$ and $H$ is an $(h+1,h+1-(n-k))$-star graph. Let
$T$ be the set of neighbors of $X$ in $S_{n,k}-X$ and $F$ be the set
of edges between $X$ and $T$. Since all the vertices with the last
$(n-1-h)$ digits $12\cdots(n-1-h)$ are in $X$, all the vertices in
$T$ are swap-neighbors of $X$ and no two vertices in $X$ share a
common swap-neighbor in $T$, that is, $|F|=|T|$.

For a vertex of $H$, since it has $h$ neighbors in $X$, it has
exactly $(n-1-h)$ neighbors in $T$. It follows that
$$
  |F|=|T|=|X|(n-1-h)=\frac{(h+1)!(n-1-h)}{(n-k)!}.
 $$

We show that $F$ is an $h$-edge-cut of $S_{n,k}$. To this end, we
only need to show that any vertex $v$ in $S_{n,k}-X$ has at least
$h$ neighbors in $S_{n,k}-F$. In fact, since $S_{n,k}$ is
$(n-1)$-regular and $v$ has at most one neighbor in $X$, $v$ has at
least $n-2\,(\geq h)$ neighbors in $S_{n,k}-X$, which implies that
$F$ is an $h$-edge-cut of $S_{n,k}$. It follows that
 $$
 \lambda_s^{(h)}(S_{n,k})\leq |F|=\frac{(h+1)!(n-1-h)}{(n-k)!}.
 $$

We now show that $T$ is an $h$-vertex-cut of $S_{n,k}$. To this end,
we only need to show that any vertex $u$ in $S_{n,k}-(X\cup T)$ has
at least $h$ neighbors within.

We claim that at most one of neighbors of $u$ is in $T$. Suppose to
the contrary that $u$ has two distinct neighbors $v$ and $w$ in $T$.
Since all vertices in $T$ are swap-neighbors of $X$, without loss of
generality, we may assume
\begin{eqnarray}
v\!\!&=&\!\!1 p_2 \ldots
p_{h+1-(n-k)}p_123\cdots (n-h-1),\qquad\qquad\label{e3.1}\\
w\!\!&=&\!\!2 p'_2 \ldots p'_{h+1-(n-k)}1p'_13\cdots
(n-h-1).\qquad\qquad\label{e3.2}
\end{eqnarray}

Since $u$ and $w$ are adjacent, their $1$-digits are different, that
is, the $1$-digit of $u$ is not $2$. If $v$ is an unswap-neighbor of
$u$, then from (\ref{e3.1}) we should have
 \begin{equation}\label{e3.3}
 u=q_1 p_2 \ldots p_{h+1-(n-k)}p_123\cdots (n-h-1),\qquad\qquad
 \end{equation}
where $q_1\in I_n\setminus\{p_1, p_2, \ldots,
p_{h+1-(n-k)},1,2,\cdots, (n-h-1)\}$. Since $p_1\ne 1, p'_1\ne 2$
and $q_1\ne 2$, Comparing (\ref{e3.2}) and (\ref{e3.3}), we can
easily find that $u$ and $w$ have different digits at least three
positions. By Definition~\ref{def2.2}, $w$ is not a neighbor of $u$,
a contradiction.

If $v$ is a swap-neighbor of $u$ then, without loss of generality,
 \begin{equation}\label{e3.4}
 u=3p_2 \ldots p_{h+1-(n-k)}p_121\cdots (n-h-1).\qquad\qquad
 \end{equation}
Comparing (\ref{e3.2}) and (\ref{e3.4}), we can also easily find
that $u$ and $w$ have different digits at least three positions, and
so $w$ is not a neighbor of $u$, a contradiction.

Since $u$ has at most one neighbor in $T$, $u$ has at least
$(n-1)-1$ neighbors in $S_{n,k}-(X\cup T)$. Since $(n-1)-1\geq h$,
$u$ has at least $h$ neighbors in $S_{n,k}-(X\cup T)$. It follows
that $T$ is an $h$-vertex-cut of $S_{n,k}$, and so
 $$
 \kappa_s^{(h)}(S_{n,k})\leq |T|=\frac{(h+1)!(n-1-h)}{(n-k)!}.
 $$
The lemma follows.
\end{pf}

\vskip6pt

By Lemma~\ref{lem3.1} and Lemma~\ref{lem3.2}, we immediately obtain
our main result.

\begin{thm}\label{thm3.3}
For $2 \leq k \leq n-1$ and $n-k\leq h \leq n-2$,
 \begin{equation}\label{e3.5}
 \kappa_s^{(h)}(S_{n,k})=\lambda_s^{(h)}(S_{n,k})=\frac{(h+1)!(n-h-1)}{(n-k)!}.
 \end{equation}
\end{thm}

\begin{rem}{\rm We would like to make some remarks on our result.

When $h=n-k$, the results in (\ref{e1.1}), (\ref{e1.2}) and
(\ref{e3.5}) are consistent, that is,
$\kappa_s^{(n-k)}(S_{n,k})=\lambda_s^{(n-k)}(S_{n,k})=(n-k+1)(k-1)$.

The alternating group network $AN_n (n\geq 3)$, proposed by
Ji~\cite{j98}, is a Cayley graph on the alternating group $A_n$ with
respect to the generating set $S=\{(1 2 3),(1 3 2),(1 2)(3 i)\ : \ 4
\leq i \leq n\}$. Cheng {\it et al.}~\cite{cqs12} proved that
$AN_n\cong S_{n,n-2}$. When $h\in\{0,1,2\}$, $\kappa_s^{(h)}(AN_n)$
and $\lambda_s^{(h)}(AN_n)$ can be obtained by (\ref{e1.1}) and
(\ref{e1.2}). Very recently, Feng {\it et al.} \cite{fhz16} have
determined $\lambda_s^{(3)}(AN_n)=12(n-4)$ for $n\geq 5$. When
$h\geq 2$, the following result is obtained from
Theorem~\ref{thm3.3} immediately.}
\end{rem}

\begin{cor}\label{cor3.5}
$\kappa_s^{(h)}(AN_n)=\lambda_s^{(h)}(AN_n)=\frac 12\,(h+1)!(n-h-1)$
for $2\leq h\leq n-2$.
\end{cor}

\section{Conclusions and Problems}

This paper considers the refined measure, $k$-super connectivity
$\kappa_s^{(h)}$ and $k$-super edge-connectivity $\lambda_s^{(h)}$
for the fault-tolerance of a network, and the $(n,k)$-star graph
$S_{n,k}$ ($2\leq k \leq n-1$), which is an attractive alternative
network to the hypercube. In early articles \cite{lx12,lx14}, we
determined $\kappa_s^{(h)}(S_{n,k})$ and $\lambda_s^{(h)}(S_{n,k})$
for $0\leqslant h\leqslant n-k$, which are two different values. In
this paper, we proved
$\kappa_s^{(h)}(S_{n,k})=\lambda_s^{(h)}(S_{n,k})=\frac{(h+1)!(n-h-1)}{(n-k)!}$
for $n-k\leq h \leq n-2$. This result implies that at least
$\frac{(h+1)!(n-h-1)}{(n-k)!}$ vertices or edges have to be removed
from an $(n,k)$-star $S_{n,k}$ to make it disconnected and no
vertices of degree less than $h$. When the $(n,k)$-star graph is
used to model the topological structure of a large-scale parallel
processing system, this result can provide a more accurate measure
for the fault tolerance of the system.

The $(n,k)$-star graph $S_{n,k}$ is vertex-transitive (see Chiang
{\it et al.}~\cite{cc95}). As we all know, $S_{n,1}\cong K_n$,
$S_{n,n-1}\cong S_n$ and $S_{n,n-2}\cong AN_n$, and so $S_{n,k}$ is
a Cayley graph for $k\in\{1,n-1,n-2\}$.

\begin{prob} Is $S_{n,k}$ a Cayley graph for
$2\leqslant k\leqslant n-3$ ?
\end{prob}

\end{document}